\documentclass[12pt]{article}

\usepackage{amsmath,amsthm,amsfonts,latexsym,amsopn,verbatim,amscd,amssymb,}
\usepackage{amssymb}
\usepackage{amsfonts}

\usepackage{graphics,color}
\usepackage{graphicx}

\theoremstyle{plain}
\newtheorem{Thm}{Theorem}[section]
\newtheorem{Lem}[Thm]{Lemma}
\newtheorem{Prop}[Thm]{Proposition}
\newtheorem{Cor}[Thm]{Corollary}
\theoremstyle{definition}

\newtheorem{Def}[Thm]{Definition}

\numberwithin{equation}{section}

\DeclareMathOperator{\idem}{Idem}
\DeclareMathOperator{\nil}{Nil}
\DeclareMathOperator{\U}{U}

\DeclareMathOperator{\en}{End}

\DeclareMathOperator{\ann}{ann}

\newcommand{\bnum}{\begin{enumerate}}
\newcommand{\enum}{\end{enumerate}}
\begin{document}

\begin{center}
\textbf{Weak Nil Clean Rings }\\
\end{center}
\begin{center}
Dhiren Kumar Basnet\\
\small{\it Department of Mathematical Sciences, Tezpur University,
 \\ Napaam, Tezpur-784028, Assam, India.\\
Email: dbasnet@tezu.ernet.in}
\end{center}
\begin{center}
Jayanta Bhattacharyya \\
\small{\it Department of Mathematical Sciences, Tezpur University,
 \\ Napaam, Tezpur-784028, Assam, India.\\
Email: jbhatta@tezu.ernet.in}
\end{center}
\noindent \textit{\small{\textbf{Abstract:}  }} We introduce the concept of a weak nil clean ring, a generalization of nil clean ring, which is nothing but a ring with unity in which every element can be expressed as sum or difference of a nilpotent and an idempotent. Further if the idempotent and nilpotent commute the ring is called weak* nil clean. We characterize all $n\in \mathbb{N}$, for which $\mathbb{Z}_n$ is weak nil clean but not nil clean. We show that if $R$ is a weak* nil clean  and $e$ is an idempotent in $R$, then the corner ring $eRe$ is also weak* nil clean.  Also we discuss $S$-weak nil clean rings and their properties, where $S$ is a set of idempotents and show that if $S=\{0, 1\}$, then a $S$-weak nil clean ring contains a unique maximal ideal. Finally  we show that weak* nil clean rings are exchange rings and strongly nil clean rings provided $2\in R$ is nilpotent in the later case. We have ended the paper with introduction of weak J-clean rings.
\bigskip

\noindent \small{\textbf{\textit{Key words:}} Nil clean ring, weak and weak* nil clean ring, exchange ring, S-weak nil clean ring, weak J- clean ring} \\
\smallskip

\noindent \small{\textbf{\textit{$2010$ Mathematics Subject Classification:}}  16N40, 16U99.} \\
\smallskip

\pagebreak

\section{Introduction}
$\mbox{\hspace{.5cm}}$ Rings $R$ are associative rings with unity unless otherwise indicated and modules  are unitary. The Jacobson radical, group of units, set of idempotents and set of nilpotent elements of  $R$ are denoted by J$(R), \U(R), \idem(R)$ and $\nil(R)$ respectively. In the  paper ``Lifting idempotents and exchange rings"\cite{wkn} Nicholson called an element $r$ in a ring $R$ clean element, if $r=e+u$ for some $e\in \idem(R)$ and $u\in \U(R),$ and a ring is clean if every element of the ring is a clean element. Similarly a nil clean ring was introduced  by Diesel \cite{ajd} in his doctoral  thesis and defined an element $r$ in a ring $R$ to be nil clean if $r=e+n$ for $e\in\idem(R)$ and $n\in \nil(R)$. A ring $R$ is nil clean if each element of $R$ is nil clean. In \cite{and} Ahn and Andreson defined a ring $R$ to be weakly clean if each element $r\in R$ can be written as $r=u+e$ or $r=u-e$ for $u\in \U(R)$ and $e\in \idem(R).$ Motivated by above concept, we observe the example $\mathbb{Z}_6=\{\overline{0},~\overline{1},~\overline{2},~\overline{3},~\overline{4},~\overline{5}\},$ here $\idem(\mathbb{Z}_6)=\{\overline{0},~\overline{1},~\overline{3},~\overline{4}\}$ and $\nil(\mathbb{Z}_6)=\{\overline{0}\}.$ So clearly $\mathbb{Z}_6$ is not a nil clean ring as $\overline{2}\mbox{ and }\overline{5}$ can not be written as sum of a idempotent and a nilpotent of $\mathbb{Z}_6.$ But we see that every elements $r\in \mathbb{Z}_6$ can be written as $r=n-e\mbox{ or }r=n+e$ for $e\in \idem(\mathbb{Z}_6)$ and $n\in \nil(\mathbb{Z}_6),$ which led us to introduce weak nil clean ring. Weak nil clean ring a ring with unity in which each element of the ring can be expressed as sum or difference of a nilpotent and an idempotent. A study on commutative weak nil clean rings have been done by Peter V. Dancheva and W. Wm. McGovernb\cite{peter} as commutative weakly nil clean rings. Here we have also given stronger version of few of its results. We have also characterized all natural number $n$, for which $\mathbb{Z}_n$ is a weak nil clean ring but not  nil clean ring.  Further we have discussed $S$-weak nil clean rings, a ring  in which each element can be expressed as sum or difference of a nilpotent and an element of $S$, where $S \subseteq \idem(R)$ and have shown that if $S=\{0, 1\}$, then a $S$-weak nil clean ring contains  unique maximal ideal. Finally we have shown that weak* nil clean rings are exchange rings and strongly nil clean rings provided $2\in R$ is nilpotent in the later case. We have ended the paper by introducing weak J-clean rings and obtain few introductory results on weak J-clean rings as a effort to answer \textbf{Problem $5$} of \cite{peter}.

\section{Weak Nil Clean Rings}
\begin{Def}An element  $r\in R$ is said to be a \textit{weak nil clean} element of the ring $R,$ if $r=n-e \mbox{ or } r=n-e,$ for some $n\in \nil(R)$, $e\in \idem(R)$ and a ring is said to be \textit{weak nil clean} ring if each of its elements is weak nil clean. Further if $r=n-e \mbox{ or }n+e$ with $ne=en$, then $r$ is called \textit{weak* nil clean}.
\end{Def}
Obviously every nil clean ring is weak nil clean, but the above example denies the converse. Also if $R$ is a weak nil clean ring or a weak* nil clean ring then for $n\geq 2,~~S=\{A = (a_{ij})\in T_n(R): a_{11}=a_{22}=\dots=a_{nn}\}$, is weak nil clean ring  which is not weak* nil clean, where $T_n(R)$ is the ring of upper triangular matrices of order $n$ over $R$. Analogue to the concept of clean and nil clean rings, it is easy to see that every weak nil clean ring is weakly clean and the converse is not true. The following theorem is easy to see.

\begin{Thm}\label{hom} Homomorphic image of a weak nil clean ring is weak nil clean.
\end{Thm}

However the converse is not true as $\mathbb{Z}_6\cong \mathbb{Z}/{\langle 6\rangle}$ is a weak nil clean ring, but $\mathbb{Z}$ is not a weak nil clean ring. A finite direct product $\prod R_\alpha$ of rings is nil clean if and only if each $R_\alpha$ is nil clean. Next result will show that similar result is not true for weak nil clean rings (following result is general form of statement (ii) of \textbf{Proposition 1.9} of \cite{peter}).

\begin{Thm}Let $\{R_\alpha\}$ be a finite collection of rings. Then the direct product $R=\prod R_\alpha$ is weak nil clean if and only if each $R_\alpha$ is weak nil clean and at most one $R_\alpha$ is not nil clean.
\end{Thm}

\noindent$Proof.$ $(\Rightarrow)$ Let $R$ be weak nil clean, then each $R_\alpha$ being homomorphic image of $R$ is weak nil clean. Suppose for some $\alpha_1$ and $\alpha_2,~~\alpha_1\neq \alpha_2,~~R_{\alpha_1}$ and $R_{\alpha_2}$ are not nil clean. Since $R_{\alpha_1}$ is not nil clean, not all elements $x \in R_{\alpha_1}$ are of the form $n - e$, where $n\in \nil(R_{\alpha_1})$ and $e\in\idem(R_{\alpha_1})$. But $R_{\alpha_1}$ is weak nil clean, so there exists $x_{\alpha_1}\in R_{\alpha_1},$ with $x_{\alpha_1}=n_{\alpha_1}+e_{\alpha_1},$ where $e_{\alpha_1}\in \idem(R_{\alpha_1})$ and $n_{\alpha_1}\in \nil (R_{\alpha_1})$, but $x_{\alpha_1}\neq n-e$ for any $n\in \nil(R_{\alpha_1})$ and $e\in\idem(R_{\alpha_1})$. Likewise there exists  $x_{\alpha_2}\in R_{\alpha_2},$ with $x_{\alpha_2}=n_{\alpha_2}-e_{\alpha_2},$ where $e_{\alpha_2}\in \idem(R_{\alpha_2})$ and $n_{\alpha_2}\in \nil (R_{\alpha_2}),$ but $x_{\alpha_2}\neq n+e$ for any $n\in \nil(R_{\alpha_2})\mbox{ and }e\in \idem(R_{\alpha_2}).$
\begin{align}\mbox{Define }x=(x_{\alpha})\in R \mbox{ by }x_{\alpha}&=x_{\alpha_i}\mbox{  if } \alpha\in \{ \alpha_1, \alpha_2\}\hspace{5cm}\notag\\\notag
          &=0 \mbox{ \hspace{.23cm}   if } \alpha \notin \{ \alpha_1, \alpha_2\}
\end{align}
Then clearly $x\neq n \pm e$ for any  $n\in \nil(R)\mbox{ and }e\in\idem(R),$ hence at most one $R_\alpha$ is not nil clean. \\
$(\Leftarrow)$ If each $R_\alpha$ is nil clean, then $R=\prod R_\alpha$ is nil clean, so weak nil clean. So assume some $R_{\alpha_0}$ is weak nil clean but not nil clean and that all other $R_\alpha$'s are nil clean. Let $x=(x_\alpha)\in R.$ In $R_{\alpha_0}$ we can write $x_{\alpha_0}=n_{\alpha_0}+e_{\alpha_0}\mbox{ or } x_{\alpha_0}= n_{\alpha_0}-e_{\alpha_0},$ where $n_{\alpha_0}\in \nil(R_{\alpha_0}),~~e_{\alpha_0}\in\idem(R_{\alpha_0}).$
If $x_{\alpha_0}=n_{\alpha_0}+e_{\alpha_0},$ for $\alpha\neq \alpha_0,$ let $x_{\alpha}=n_{\alpha}+e_{\alpha}$ and if $x_{\alpha_0}=n_{\alpha_0}-e_{\alpha_0},$ for $\alpha \neq \alpha_0,$ let $x_{\alpha}=n_{\alpha}-e_{\alpha}$ then $n=(n_\alpha)\in\nil(R) \mbox{ and } e=(e_\alpha)\in\idem(R) \mbox{ and } x=n+e \mbox{ or } x=n-e$ respectively, hence $R$ is weak nil clean. \hfill $\square$

\begin{Prop}\label{12}Let $R$ be a weak nil clean ring, then $J(R)\subseteq \nil(R).$
\end{Prop}

\noindent$Proof.$ Let $x\in J(R).$ Then $x=n-e\mbox{ or }x=n+e,$ where $n\in \nil(R) \mbox{ and } e\in\idem(R)$. If $x=n-e$ then there exists a $k\in \mathbb{N}$ such that $(x+e)^k=0$, which gives $ e \in J(R) \cap\idem(R)$, hence $e = 0$ i.e., $ x = n\in\nil(R).$ Similarly for $x=n+e$, we get $x=n\in\nil(R)$. Thus $ J(R)\subseteq \nil(R).$ \hfill $\square$\\


\begin{Prop}If a commutative ring $R$ is weak nil clean, $R/\nil(R)$ is weak nil clean and converse holds if idempotents can be lifted modulo $\nil(R).$
\end{Prop}

\noindent$Proof.~~(\Rightarrow)$ Follows from Theorem \eqref{hom}.\\
$(\Leftarrow)$ Let $x\in R$. Since $R/\nil(R)$ is weak nil clean, so
$x+\nil(R)=y+\nil(R) \mbox{ or } (-y)+\nil(R),$ where
$y^2-y\in \nil(R$) ( as $R/\nil(R)$ is a reduced ring). Since idempotents of $R$ lift modulo $\nil(R)$, so there exist $e\in\idem(R)$ such that $y-e\in\nil(R)$, which implies $x-e\in\nil(R) \mbox{ or }x+e\in\nil(R)$ i.e.,
$x-e=n \mbox{ or }x+e=m \mbox{ for some } m,n\in \nil(R),$ which proves the result. \hfill $\square$\\

\indent For more examples of weak nil clean rings, we consider the method of idealization. Let $R$ be a commutative ring and $M$  a left $R-$module. The idealization of $R$ and $M$ is the ring $R(M)=R\oplus M$ with product defined as $(r,m)(r^\prime,m^\prime)=(rr^\prime,rm^\prime+r^\prime m)$ and sum as $(r,m)+(r^\prime,m^\prime)=(r+r^\prime,m+m^\prime),$ for $(r,m),~~(r^\prime,m^\prime)\in R(M).$

\begin{Thm}Let $R$ be a ring and $M$ be a left $R$-module. Then $R$ is weak nil clean if and only if $R(M)$ is weak nil clean.
\end{Thm}

\noindent$Proof.~~(\Leftarrow)$ Note that $R\approx R(M)/(0\oplus M)$ is homomorphic image of $R(M)$. Hence by Theorem \eqref{hom},  $R$ is weak nil clean ring.\\
$(\Rightarrow)$ Let $R$ be weak nil clean ring and $(r,m)\in R\oplus M,$ where $r\in R$ and $m\in M,$ we have $r=n+e\mbox{ or } n-e$ for $n\in\nil(R)\mbox{ and }e\in\idem(R)$, then
$(r,m)=(n+e,m)\mbox{ or }(n-e,m)=(n,m)+(e,0)\mbox{ or }(n,m)-(e,0)$ is weak nil clean expression of $(r,m)$, where $(n,m)\in\nil(R)\mbox{ and }(e,0)\in\idem(R),$ hence $R(M)=R\oplus M$ is weak nil clean. \hfill $\square$\\

Now we try to characterize all $n$ for which $ \mathbb{Z}_n$ is weak nil clean but not nil clean. We recall that, $\idem{(\mathbb{Z}_{p^k})}=\{0,1\}$, for any prime $p\in \mathbb{N}$ and $k\in \mathbb{N}$.
\begin{Lem} \label{11}
$\mathbb{Z}_{3^k},~~k\in \mathbb{N}$ is weak nil clean but not nil clean.
\end{Lem}
\noindent$Proof$. The proof follows from the fact that $\idem{(\mathbb{Z}_{3^k})}=\{0,1\}$ and $\nil{(\mathbb{Z}_{3^k})}=\{0,3, 6, . . ., 3(3^{k-1}-1)\}$. 
\begin{Lem}  $\mathbb{Z}_{p^k}, k \in \mathbb{N}$ is weak nil clean but not nil clean, where $p$ is prime iff $p=3$.
\end{Lem}
\noindent$Proof$. $(\Leftarrow)$ It follows from Lemma \ref{11}\\
$(\Rightarrow)$ We know that  $\mathbb{Z}_{2^k}$ is nil clean $\forall k\in \mathbb{N}$ and  $\mathbb{Z}_{3^k}$ is weak nil clean $\forall k\in \mathbb{N}$ but not nil clean. Now consider $p>3$, then we have $\idem{(\mathbb{Z}_{p^k})}=\{0,1\}$ and $\nil{(\mathbb{Z}_{p^k})}=\{0,p,2p,\dotsc,(p^{k-1}-1)p\}$.
So if we consider the sum or difference of nilpotents and idempotents of $\mathbb{Z}_{p^k}$ respectively, then at most $4p^{k-1}$ elements can be obtained, but $p>4$, so $p^k>4p^{k-1}$. Hence all elements of $\mathbb{Z}_{p^k}$ can not be written as a sum or difference of nilpotent and idempotent of $\mathbb{Z}_{p^k}$ respectively. So $p=3.\hfill\square$
\begin{Thm} The only $n$ for which $\mathbb{Z}_n$ is weak nil clean but not nil clean is of the form  $2^r3^t$, where $t \in \mathbb{N}, r \in \mathbb{N}\cup \{0\}$. 
\end{Thm}
\noindent$Proof.$ We have already seen that $\mathbb{Z}_{3^t}$ is weak nil clean but not nil clean. Next let $n=p^{\alpha_1}_1p^{\alpha_2}_2\dotsm p^{\alpha_k}_k$ with $\alpha_i\in \mathbb{N},~1\leq i\leq k$ and $p_i$'s are distinct primes such that $p_1\leq p_2\leq . . . \leq p_n$.  If $k>2$, then there exists some $i$ with $1\leq i\leq k$ such that $p_i>3$. Then $\mathbb{Z}_{p_i^{\alpha_i}}$ is not weak nil clean. Hence $\mathbb{Z}_n$ can not be weak nil clean as $\mathbb{Z}_n=\mathbb{Z}_{p_1^{\alpha_1}}\oplus \mathbb{Z}_{p_2^{\alpha_2}} \oplus\dotsb \oplus\mathbb{Z}_{p_k^{\alpha_k}}$. So $k\leq 2$ and $p_i\leq 3$ i.e., $n=p_1^{\alpha_1}p_2^{\alpha_2}$. If $k=1$, then $p_1$ must be $3$ as $\mathbb{Z}_{2^r}$ is nil clean. Again if $k=2$, then since $p_i$'s are distinct so $p_1=2$ and $p_2=3$. Also if $n=2^{\alpha_1}3^{\alpha_2}$, then $\mathbb{Z}_n=\mathbb{Z}_{2^{\alpha_1}}\oplus \mathbb{Z}_{3^{\alpha_2}}$. Since $\mathbb{Z}_{2^{\alpha_1}}$ is nil clean and $\mathbb{Z}_{3^{\alpha_2}}$ is weak nil clean but not nil clean, so $\mathbb{Z}_n$ is weak nil clean but not nil clean. This completes the proof.$\hfill\square$\\



 The polynomial ring $R[x]$ over a weak nil clean ring is not necessarily weak nil clean. In fact if $R$ is commutative the $R[x]$ is never weak nil clean. For then $x \in R[x]$ is of the form $\sum_ia_ix^i-e$ or $\sum_ia_ix^i+e$, where $a_i \in \nil(R), e \in \idem(R)$, giving $a_0-e=0$ or $a_0+e=0$, which is absurd.\\

\noindent However if $R$ is weak nil clean and $\sigma : R \rightarrow R$ is a ring endomorphism then for any $n$, the quotient $S = R[x;\sigma]/<x^n>$, where $R[x;\sigma]$ is the Hilbert twist, is a weak nil clean ring. Indeed if $f=a_0+a_1x+a_2x^2+ . . . +a_{n-1}x^{n-1} \in S$ and $a_0 = n+ e$ or $a_0= n - e$, where $n \in \nil(R), e \in \idem(R)$, then $f=(f-e)+e$ or $f = (f+e)-e$ is a weak nil clean decomposition of $f$ in $S$.\\

In order to show that, weak* nil cleanness penetrates to corner, we need following lemmas.

\begin{Lem} \label{1} Let $R$ be a ring and $x=n+e$ or $n-e$ be a weak* nil clean decomposition of $x\in R$ with $n\in \nil(R)$ and $e\in \idem(R),$ then
$\ann _l(x)\subseteq \ann _l(e)$ and $\ann _r(x)\subseteq \ann _r(e)$, where $\ann _l(a)$ and $\ann _r(a)$ denote the left and right annihilator of an element $a$ in $R$ respectively.
\end{Lem}

\noindent$Proof.$ Let $r\in\ann_l(x)$ then $rx=0.$ If $x=e+n$ then $rn+re=0$, and so $ rne+re=0$ i.e., $ re(n+1)=0$  implying $re=0$ and thus $ r\in \ann_l(e)$.\\
Again  if $x=n-e$, then $rn-re=0$ and so $rne-re=0$ i.e., $re(n-1)=0$ implying $re=0$ and thus $r\in \ann_l(e)$.
Hence $\ann _l(x)\subseteq \ann _l(e)$. Similarly the other part i.e., $\ann _r(x)\subseteq \ann _r(e)$. \hfill $\square$

\begin{Lem}\label{2} Let $R$ be a ring and $x=n+e$ or $n-e$ be a weak* nil clean decomposition of $x\in R$ with $n\in \nil(R)$ and $e\in \idem(R),$ then
$\ann_l(x)\subseteq R(1-e)$ and $\ann_r(x)\subseteq (1-e)R.$
\end{Lem}
\noindent$Proof.$ Straightforward.

\begin{Thm}Let $R$ be a ring and $f\in Idem(R),$ then $x\in fRf$ is weak* nil clean in $R$ if and only if $x$ is weak* nil clean in $fRf.$
\end{Thm}
\noindent$Proof.(\Leftarrow)$
If $x\in fRf$ is weak* nil clean in $fRf$, then by the same weak* nil clean decomposition, $x$ is weak* nil clean in R.\\
$(\Rightarrow)$ Let $x$ is weak* nil clean in $R,$ so $x=n+e \mbox{ or }n-e$ for some $n\in\nil(R)\mbox{ and }e\in\idem(R)$ with $ne=en$.  First let $x=n+e$, since $x\in fRf,$ so $(1-f)\in \ann_l(x)\cap \ann_r(x)\subseteq R(1-e)\cap (1-e)R = (1-e)R(1-e)$ [ by Lemma \ref{2}]. Thus we have $(1-f)e=0=e(1-f)$ giving $fe=e=ef$ and consequently $fef\in\idem(fRf).$ Also $xf=fx$, therefore we have $nf=fn$ i.e., $fnf\in \nil(fRf)$. Hence $ x= fnf+fef$. Similarly if $x=n-e$ then $x=fnf-fef.$ Hence $x$ is weak* nil clean in $fRf$.\hfill $\square$\\

The following is an immediate consequence of this result.\\

\begin{Cor}Let $R$ be weak* nil clean ring and $e\in\idem(R)$, then the corner ring $eRe$ is also weak* nil clean.
\end{Cor}

\section{S-Weak Nil Clean Rings}
S-weak nil clean ring is a generalization of weak nil clean ring, which is defined as follows:
\begin{Def} let $S$ be a non-empty set of idempotents of $R$, then $R$ is called $S-weak~~nil~~clean$ if each $r\in R$ can be written as $r=n+e\mbox{ or }n-e,$ where $n\in\nil(R)\mbox{ and }e\in S.$
\end{Def}
\begin{Prop}
Let $R$ be $\{0,1\}$-weak nil clean ring, then $R$ has exactly one maximal ideal.
\end{Prop}
\noindent$Proof$.
 Since $R$ is $\{0,1\}$-weak nil clean ring so $R=\U(R) \bigcup \nil(R)$ and $\U(R)=(1+\nil(R))\bigcup (-1+\nil(R))$. It follows that for any $x\in \nil(R)$ and any $r\in R$, we have $xr, rx \in \nil(R).$
Next if possible let $n_1 -n_2=u$, where $n_1,n_2\in\nil(R) \mbox{ and }u\in\U(R)$. Then $u^{-1}n_1-u^{-1}n_2=1$ i.e., $ n_3=1+n_4$, where $u^{-1}n_1=n_3\in\nil(R)$ and $u^{-1}n_2=n_4\in \nil{(R)}$, which is a contradiction. Thus $n_1 -n_2\in\nil(R)$, for any $n_1,~~n_2\in\nil(R)$ implying that $\nil(R)$ is an ideal. Hence by proposition \ref{12} $J(R) = \nil(R)$. This completes the proof.\hfill $\square$\\

\indent From above theorem it is clear that $\{0,1\}-$ nil clean rings are local rings. Converse is not true.

\begin{Thm}If a ring $R$ is S-weak* nil clean for $S\subseteq \idem(R)$ then $S = \idem(R)$.
\end{Thm}

\noindent$Proof.$ Let $e^\prime\in\idem(R)$, then $-e^\prime\in R$. Since $R$ is $S$-weak* nil clean, so $-e^\prime=n+e\mbox{ or }-e^\prime=n-e$ for some $n\in\nil(R),\mbox{ and }e\in S$, with $ne=en$.
If $-e^\prime=n+e$, then $1-e^\prime =1+n+e$ i.e., $(1+n+e)^2=1+n+e$, which gives $1+n^2+e+2n+2e+2ne=1+n+e$ i.e., $n^2+n+2e(1+n)=0$, implies $(n+2e)(1+n)=0$. But $1+n\in U(R)$ , so $n=-2e$, giving $-e^\prime=n+e=-2e+e=-e$. Thus $e^\prime=e\in S$. \\
Again if $-e^\prime=n-e$, then $(-e^\prime)^2={e^\prime}^2 =e^\prime$ i.e., $(n-e)^2=-n+e$, which gives $n^2-2ne+e=-n+e$ i.e., $n^2+n(1-2e)=0$, implies $n\{n+(1-2e)\}=0$. But $n+(1-2e)\in U(R)$ , so $n=0$ i.e., $e^\prime=e\in S$. Hence $\idem(R)=S.\hfill\square$ \\

But in case of weak clean ring it is possible that $R$ is $S-$weak clean and $S\subsetneq \idem(R)$ \cite{and}.

\section{More result on  Weak Nil Clean Rings}

\indent It is well known that $\mathbb{Z}_3$ is clean, so upper triangular matrix ring $\mathbb{T}_2(\mathbb{Z}_3)$ is clean and hence exchange, but  $\mathbb{T}_2(\mathbb{Z}_3)$ is not weak nil clean ring, so in general, exchange rings are not weak nil clean rings. But one can see that weak* nil clean rings are exchange.

\begin{Thm}Let $R$ be a weak* nil clean ring, then $R$ is a exchange ring.
\end{Thm}

\noindent$Proof.$ Let $R$ be a weak* nil clean ring and $x\in R,$ then $x=n+e\mbox{ or } x=n-e,$ where $n\in\nil(R)\mbox{ and }e\in \idem(R).$\\
If $x=n-e$, then
\begin{align}(1-n)[x-(1-n)^{-1}e(1-n)]&=(1-n)[(n-e)-(1-n)^{-1}e(1-n)],\notag\\
                                      &=n-e-n^2+ne-e+en,\notag\\
                                      &=x-(n-e)^2=x-x^2,\notag\\
\mbox{implies} \hspace{1.5cm}[x-(1-n)^{-1}e(1-n)]&=(1-n)^{-1}(x-x^2).\notag
\end{align}
Similarly if $x=n+e$, we have $x-e=u^{-1}(x^2-x)\mbox{ for } u=(2e-1)+n \in \U(R)$. Then by  condition $(1)$ of Proposition $1.1$ of \cite{wkn}, $R$ is exchange. \hfill $\square$\\

\indent Finally we take the question `` under what condition a weak* nil clean ring is strongly nil clean ring?" To answer this question we need the following  Lemma.

\begin{Lem}\label{wnc}Let $R$ be a ring and $M_R$  a right $R$-module. If an endomorphism $\phi\in\en(M_R)$ is sum or difference of a nilpotent $n$ and an idempotent $e$, which commutes with $2\in Nil(R)$ then there exists a direct sum decomposition $M=A\oplus B$ such that $\phi|_A$ is an element of $\en(A)$ which is nilpotent and $(1-\phi)|_B$ is an element of $\en(B)$ which is nilpotent.
\end{Lem}

\noindent$Proof.$ Suppose $\phi=a-e,$ where $e\in\idem(\en(M_R))\mbox{ and } a\in\nil(\en(M_R))$ and suppose $ea=ae.$ We define decomposition $M=A\oplus B,$ by setting\\
$A=(1-e)M$ and $B=eM.$ Then $A$ and $B$ are $\phi-$invariant.\\
Now $\phi|_A=(a-e)|_A=a|_A-e|_A=a|_A$ and so $\phi|_A$ is nilpotent.\\
And $(1-\phi)|_B=(1-(a-e))|_B=(1-a+e)|_B=(2-a - (1-e))|_B =(2-a)|_B$ is nilpotent as $2$ is nilpotent. \\
Again, if $\phi=a+e,$ where $e\in\idem(\en(M_R))\mbox{ and } a\in\nil(\en(M_R))$, then by Definition $1.2.8$ and Lemma $1.2.3$ of \cite{ajd} such a decomposition exists. $\hfill \square$ \\

\indent Now we can state following theorem.

\begin{Thm}A ring $R$ is strongly nil clean if and only if $R$ is weak* nil clean with $2\in \nil(R)$.
\end{Thm}
\noindent$Proof.(\Rightarrow)$ Clear from the definition of weak* nil clean ring.\\
$(\Leftarrow)$ The result follows from Lemma $1.2.6$ of \cite{ajd} and Lemma \eqref{wnc}.

\begin{Cor}A weak* nil clean ring $R$ with $2\in \nil(R),$ is strongly $\pi-$regular.
\end{Cor}

\section{Weak J-clean ring}
In this section we have defined weak J-clean element of a ring, as a generalization of J-clean rings by Chen \cite{hc3}.
\begin{Def}
An element $a$ in a ring $R$ is said to be weak J-clean if $a$ can be written as $a=j+e$ or $a=j-e$ for some $j\in J(R)$ and $e\in \idem(R)$. Moreover if
$ae=ea$ we say a to be weak* J-clean.
\end{Def}
Following are some of the preliminary result we got related to weak J-clean rings.
\begin{Lem}
Every weak* J-clean element in a ring is strongly clean.
\end{Lem}
\noindent$Proof.$ Let $a\in R$ be a ring element, $e\in \idem{R}$ and $w\in J{R}$. If $a=w+e$ we have $a=(1-e)+(2e-1+w)$, else if $a=w-e$ we have $a=(1-e)-(1-w)$.\hfill $\square$

\begin{Lem}Let $R$ be a ring  and $a=w+e$ or $a=w-e$ be weak* J-clean decomposition of $a$ in $R$, where $e\in \idem(R)$ and $w\in J(R)$. Then $\ann_l(a)\subseteq \ann_l(e)$ and $\ann_r(a)\subseteq \ann_r(e)$.
\end{Lem}
\noindent$Proof.$ Let $r\in \ann_l(a)$, then $ra=0$ consider $a=w+e$ then $re=-rw$ hence $re=-rwe=-rew$.It follows that $re=0$, thus $r\in\ann_l(e)$. Similarly $\ann_r(a)\subseteq \ann_r(e)$ can be shown. \hfill $\square$

\begin{Thm}Let $R$ be a ring and let $f\in R$ be an idempotent. Then
$a\in fRf$ is weak* J-clean in $R$ if and only if $a$ is weak* J-clean in $fRf$.
\end{Thm}
\noindent$Proof.$ Let $a\in fRf$ and $a=w+e\mbox{ or }w-e$ for $w\in J(R)$ and $e\in \idem(R)$. We begin by showing that $e\in \idem(fRf)$, then by using above weak* J-clean expression of $a$ in R, it is easy to deduce that $w\in J(fRf)$, implies that above weak* J-clean expression of $a$ is also the a weak* J-clean expression of $a$ in $fRf$. To show $e\in \idem(R)$ observe that  $1-f\in \ann_l(a)\cap \ann_r(a) \subseteq \ann_l(e)\cap \ann_r(e)$, implies $ef=e=fe$, hence $e\in \idem(R)$. Other part of the theorem follows trivially.\hfill $\square$

\begin{Cor}
Let $R$ be weak* J-clean ring, $e\in R$ be an idempotent then so is $eRe$.
\end{Cor}
Before proceeding further we have generalized one popular concept  lifting of idempotent modulo ideal $I$ of a ring $R$.
\begin{Def}
Let $I$ be an ideal of $R$. We say idempotents lift weakly modulo $I$, if for each idempotent $\overline{e}\in R/I$, there exists an idempotent $e\in R$ such that $e-f\in I\mbox{ or } e+f\in I$.
\end{Def}
\begin{Thm}
$R$ be a ring such that, $R/J(R)$ is boolean and each idempotent lifts weakly modulo $J(R)$ then $R$ is weak J-clean.
\end{Thm}
\noindent$Proof.$ For $a\in R$, $\overline{a}\in J(R)$ is an idempotent. By assumption we can find an idempotent $e\in R$, such that $a-e\in J(R)
\mbox{ or } a+e\in J(R)$. In both the cases we get a weak J-clean expression for a in $R$, Hence $R$ is weak J-clean. \hfill $\square$

\begin{Thm}
Let $R$ be a weak* J-clean ring, such that $R/J(R)$ is boolean then $R$ is J-clean.
\end{Thm}
\noindent$Proof.$ For $a\in R$ we have at least one $e\in \idem(R)$ such that $u=a+e \in U(R)$ or $u=a-e\in U(R)$. As $\overline{u}^2=\overline{u}$, we deduce that $u\in 1+J(R)$. Thus $a=(1-e)+(2e+(x-e-1)) \mbox{ or } a=(1-e) +(2e+(x+e-1))$ is a strongly J-clean decomposition of $a$.\hfill $\square$

\end{document}